\documentclass[a4paper,11pt]{amsart}
\usepackage[matrix,arrow,tips,curve]{xy}
\usepackage[english]{babel}
\usepackage{amsmath,stackrel}
\usepackage{amssymb}
\usepackage{tikz}
\usepackage{mathrsfs}
\usepackage{enumerate}
\usepackage{graphicx}
\usepackage{hyperref}
\usetikzlibrary{trees}
\usetikzlibrary{arrows}
\usepackage{xcolor}
\makeindex

\newtheorem{thm}{Theorem}[section]
\newtheorem{Lemma}[thm]{Lemma}
\newtheorem{Proposition}[thm]{Proposition}

\newtheorem{Remark}[thm]{Remark}

\newtheorem{say}[thm]{}

\def\defi[#1]{{\bf{#1}}}

\newcommand{\cA}{{\bf a}}

\renewcommand{\P}{\mathbb{P}}

\DeclareMathOperator{\tr}{tr}

\DeclareMathOperator{\C}{\mathbb{C}}

\begin{document}
\title[A Torelli theorem for parabolic vector bundles over an elliptic curve]{A Torelli theorem for moduli spaces of parabolic vector bundles over an elliptic curve}
\author[Thiago Fassarella]{Thiago Fassarella}
\address{\sc Thiago Fassarella\\
Universidade Federal Fluminense, Rua Alexandre Moura 8 - S\~ao Domingos, 24210-200 Niter\'oi, Rio de Janeiro, Brazil}
\email{tfassarella@id.uff.br}
\author[Luana  Justo]{Luana  Justo}
\address{\sc Luana Justo\\
Universidade Federal Fluminense, Rua Alexandre Moura 8 - S\~ao Domingos, 24210-200 Niter\'oi, Rio de Janeiro, Brazil
}
\email{ljusto@ifes.edu.br}

\keywords{Moduli of parabolic bundles, Higgs bundles, Torelli theorem. Luana Justo was partially supported by CAPES. The authors also thank CAPES-COFECUB Ma932/19}

\begin{abstract}
Let $C$ be an elliptic curve, $w\in C$,  and let $S\subset C$ be a finite subset of cardinality at least $3$.  We prove a  Torelli type theorem for the moduli space of rank two parabolic vector bundles with determinant line bundle $\mathcal O_C(w)$ over $(C,S)$ which are semistable with respect to a weight vector $\big(\frac{1}{2}, \dots, \frac{1}{2}\big)$. 
\end{abstract}

\maketitle

%
%

\section{Introduction}
Let $C$ be a smooth complex curve of genus $g\ge 2$ and fix $w\in C$. Let $\mathcal M$ be the corresponding moduli space of semistable rank two vector bundles having $\mathcal O_C(w)$ as determinant line bundle. A classical  Torelli type  theorem  of D. Mumford and P. Newstead \cite{MN68} says that the isomorphism class of $\mathcal M$ determines the isomorphism class of $C$. This result has been extended first in \cite{Ty70,NR75,KP95} to higher rank and later to the parabolic context, which we now describe.

We now assume $g\ge 0$. Let $S\subset C$ be a finite subset of cardinality $n\ge 1$. 
Let $\mathcal M_{\cA}$ be the moduli space of rank two parabolic vector bundles on $(C,S)$ with fixed determinant line bundle $\mathcal O_C(w)$, and which are $\mu_{\cA}$-semistable, see Section \ref{quasiparvec}.   The subscript $\cA$ refers to a particular choice of a  weight vector $\cA = (a_1,\dots, a_n)$ of real numbers, $0\le a_i\le 1$, which gives the slope-stability condition. The moduli space associated to the central weight $\cA_F = \big(\frac{1}{2}, \dots, \frac{1}{2}\big)$ is particularly  interesting, for instance when $g=0$ and $n\ge 5$ it is a Fano variety that is smooth if $n$ is odd and has isolated singularities if $n$ is even, see \cite{Mu05, Casagrande, AM16, AFKM19}. In this context of parabolic bundles, a  Torelli type theorem has been obtained in \cite{BHK10} for $g=0$ and weight vector $\cA_F$, in \cite{BBR01} for $ g\ge 2$ and small system of weights  and \cite{AG19} deals with the case $g\ge 4$  and arbitrary rank.  We focus on the case $g=1$ (see Theorem \ref{main inside}):

%
%
%

\begin{thm}\label{main}
Let $C$ and $C'$ be two elliptic curves. Let $S\subset C$ and  $S'\subset C'$ be subsets of cardinality $n\ge 3$. Given $w\in C$ and $w'\in C'$, consider the two moduli spaces $\mathcal M_{\cA_F}$ and $\mathcal M'_{\cA_F}$ of $\mu_{\cA_F}$-semistable parabolic vector bundles,  over $(C,S)$ and $(C',S')$, with fixed determinant line bundle $\mathcal O_C(w)$ and $\mathcal O_{C'}(w')$,  respectively. If these moduli spaces are isomorphic, then there is an isomorphism $C\longrightarrow C'$ sending $S$ to $S'$. 
\end{thm}

The assumption on the number $n$ of parabolic points is necessary, when $n=2$ the corresponding moduli space $\mathcal M_{\cA_F}$  is isomorphic to $\mathbb P^1\times\mathbb P^1$, see \cite{N19}. 

The main ingredient in the proof of Theorem \ref{main} is a classification, up to elementary transformations,  of indecomposable parabolic vector bundles on $(C,S)$ which are not $\mu_{\cA_F}$-stable, see Proposition \ref{classi_ind_uns}. As consequence of this proposition we show that the complement of  $T^*\mathcal M_{\cA_F}^s$ in $\mathcal H_{\cA_F}^s$ has codimension at least $\frac{n}{2}$, where $\mathcal H_{\cA_F}^s$ denotes the moduli space of $\mu_{\cA_F}$-stable parabolic Higgs bundles over $(C,S)$ which are trace-free, see  Proposition \ref{compTM}, and  the rest of the proof of Theorem \ref{main} follows the approach of \cite{HR04, BGM13, AG19}. 

\medskip

\medskip

\noindent {\bf Notation and conventions.} Throughout the paper we work over the field $\C$ of complex numbers. When $L$ is a line bundle over a variety $X$, we  often write $\Gamma(L)$ for ${\rm H}^0(X,L)$. If $S$ is a divisor on $X$ we write $L(S)$ for $L\otimes \mathcal O_X(S)$.

\medskip


%

%
%

\section{Rank two parabolic vector bundles}\label{quasiparvec}

Let $C$ be a smooth irreducible complex curve of genus $g\ge 0$. Fix  $p_1, \dots, p_n\in C$ distinct points and denote by $S=p_1+\cdots +p_n$ the effective reduced divisor determined by them. By abuse of notation we often write $S$ for $\{p_1,\dots, p_n\}$.   Let $\omega_C$ denote the canonical sheaf of $C$. 

A \textit{quasiparabolic vector bundle} $E_{\bf v} = (E, {\bf v})$, ${\bf v} = \{V_{i}\}$,  of rank two  on $\big(C, S\big)$ consists  of
a holomorphic vector bundle $E$ of rank two on $C$ and for each  $i = 1,\dots,n$, a $1$-dimensional linear subspace $V_{i} \subset E_{p_{i}}$. The points $p_i$'s are called parabolic points, and  the subspace $V_{i} \subset E_{p_{i}}$ is the parabolic direction of $E$ at $p_i$.

%

%
%

Fix a weight  vector $\cA = (a_{1}, \dots, a_{n})$ of real numbers $0 \leq a_{i} \leq 1$.
The  \textit{parabolic slope}  of $(E, {\bf v})$ with respect to $\cA$ is 
$$
\mu_{\cA}(E) = \frac{\deg E + \sum_{i=1}^{n}a_{i}}{2}
$$ 
where $\deg E = \deg (\det E)$.  Let $L \subset E$ be a line subbundle. For each  $i = 1,\dots,n$, set 
$$a_i(L,E)\ \ = \  \left\{ 
\begin{aligned}
& a_i \ & \text{ if } L_{p_{i}} = V_{i},\\
&0 \ &  \text{ if } L_{p_{i}} \neq V_{i}.
\end{aligned}
\right.$$
The  \textit{parabolic slope}  of $L \subset E$ with respect to $\cA$ is 
$$
\mu_{\cA}(L,E) =  \deg(L)+\sum_{i=1}^{n}a_i(L,E).
$$

A quasiparabolic vector bundle $(E,{\bf v})$  is $\mu_{\cA}$-\textit{semistable} (respectively $\mu_{\cA}$-\textit{stable}) if for every  line subbundle $L \subset E$ we have $\mu_{\cA}(L,E) \leq  \mu_{\cA}(E)$ (respectively $\mu_{\cA}(L,E) < \mu_{\cA}(E)$). A \textit{parabolic vector bundle} is a quasiparabolic vector bundle together with a  weight vector $\cA$. Two parabolic vector bundles are called  ${\rm S}$-equivalent if their associated graded bundles are isomorphic, see for example \cite{MS80}. 

For each fixed $d\in\mathbb Z$, there is a moduli space $\mathcal{M}_{\cA}(d)$ parametrizing ${\rm S}$-equivalence classes of rank two quasiparabolic vector bundles $E_{\bf v} = (E, {\bf v})$ on $\big(C, S\big)$, with $\deg E = d$,  which are $\mu_{\cA}$-semistable.  If $L$ is a  line subbundle with $\deg L=d$, we denote by $\mathcal{M}_{\cA}(L)$ the subvariety of $\mathcal{M}_{\cA}(d)$ given by those parabolic vector bundles with $\det E = L$. 

The moduli space $\mathcal{M}_{\cA}(d)$ is a normal projective variety of dimension $n-3+4g$, if it is not empty, see \cite{MS80,Yo95,Bho96}, whereas $\mathcal{M}_{\cA}(L)$ has dimension $n-3+3g$. By twisting vector bundles with a fixed line bundle $L_0$, we see that $\mathcal{M}_{\cA}(L) \cong \mathcal{M}_{\cA}(L\otimes L_0^2)$. 

From now on we assume that $C$ is an elliptic curve.  Given $w,w'\in C$ we can see that 
\[
\mathcal{M}_{\cA}(\mathcal O_C(w))\simeq \mathcal{M}_{\cA}(\mathcal O_C(w')).
\]
Indeed, the isomorphism is obtained by twisting vector bundles with a fixed line bundle  $L_0$, where $L_0$ is a square root of $\mathcal O_C(w'-w)$.   We write simply $\mathcal{M}_{\cA}$ for the corresponding moduli space   
\[
\mathcal{M}_{\cA} = \mathcal{M}_{\cA}(\mathcal O_C(w)).
\]
Let $\mathcal{M}_{\cA}^s$ be the Zariski open subset parametrizing stable parabolic vector bundles. The central weight
\[
\cA_F = \left(\frac{1}{2}, \dots, \frac{1}{2}\right) 
\] 
plays an important role in this work.

There is a correspondence between quasiparabolic vector bundles, called elementary transformation,  which we now describe. Given a subset $I\subset \{1,\dots, n\}$ of even cardinality we consider the following exact sequence of sheaves
$$
0\ \rightarrow\ E' \ \stackrel{\alpha}{\to} \  E\ \stackrel{\beta}{\rightarrow}\ \bigoplus_{i\in I} E/V_i \ \rightarrow\ 0 \ 
$$
where $E/V_i$ intends to be a skyscraper sheaf determined by $E_{p_i}/V_i$,  i.e., for an open subset $U$ of $C$ we have  $ (E/V_i) (U) = E_{p_i}/V_i$ if $p_i\in U$ and $\{0\}$ otherwise. The map $\beta$ sends $s$ to $\oplus_{i\in I} s(p_i)$. If $E$ is locally generated by $e_1,e_2$ as $\mathcal O_C$-module  near $p_i$ with $e_1(p_i)\in V_i$, then $E'$ is locally generated by $e_1,e_2'$, with $e_2' = xe_2$, where $x$ is a local parameter. In particular $E'$ is locally free of rank two. We view $E'$ as a quasiparabolic vector bundle $(E',{\bf v'})$ of rank two over $(C,S)$ putting $V_i' := {\rm ker}\alpha_{p_i}$.  Notice that we have the following equality
\[
\det E' = \det E \otimes \mathcal O_{C}\big(-D\big).
\]
where $D=\sum_{i\in I}p_i$. The stability condition is preserved after an appropriate 
modification of weights,  if $(E,{\bf v})$ is $\mu_{\cA}$-semistable then $(E',{\bf v'})$ is $\mu_{\cA'}$-semistable with $\cA'_i=1-a_i$ if $i\in I$ and $\cA'_i = a_i$ otherwise. 
 In particular, since $D$ has even degree, by  choosing  a square root $L_0$ of $\mathcal O_C(D)$ we obtain an isomorphism between moduli spaces $elem_{I,L_0}:\mathcal{M}_{\cA} \to \mathcal{M}_{\cA'}$
\[
elem_{I,L_0}: (E,{\bf v}) \to (E', {\bf v}') \otimes L_0  
\]
When $\cA = \cA_F$ it turns out to be an automorphism of $\mathcal{M}_{\cA_F} $, we call it an {\it elementary transformation} over $I$. We denote by ${\bf El}\subset Aut(\mathcal{M}_{\cA_F} )$ the group of all elementary transformations. 

When $I$ has odd cardinality the parity of the degree of $E$ is modified, for instance by performing an elementary transformation centered at one parabolic point $p_i$, we see that $\mathcal{M}_{\cA}(L) \cong \mathcal M_{\cA^{i}}(L\otimes \mathcal O_C(-p_i))$, where  
$$
\cA^{i}=(a_1, \dots, 1-a_i, \dots, a_n).
$$

\section{Rank two parabolic Higgs bundles}

\begin{say}{\bf Parabolic Higgs bundles.}\label{higgs}\rm
\;Let $E_{\bf v}$  be a quasiparabolic vector bundle on $(C,S)$. An endomorphism $f : E \longrightarrow E$ is called \textit{parabolic} if $f(V_{i}) \subseteq V_i$ for every $i = 1,\dots,n$. We denote by ${\it \mathcal {E}nd}(E_{\bf v})$ the sheaf of parabolic endomorphisms.  A parabolic endomorphism is called \textit{strongly parabolic} if $f(E_{p_i}) \subseteq V_i$ and $f(V_i)=0$ for every $i = 1,\dots,n$. The sheaf of strongly parabolic endomorphisms of $E_{\bf v}$ will be denoted by $\mathcal{SE}nd(E_{\bf v})$. We denote by ${\it \mathcal {E}nd}_0(E_{\bf v})$ and $\mathcal{SE}nd_0(E_{\bf v})$ the sheaves of parabolic and strongly parabolic endomorphisms of vanishing trace.

A \textit{traceless Higgs field} on $E_{\bf v}$ is a section
$$
\theta \in \Gamma(\mathcal{SE}nd_0(E_{\bf v})\otimes \omega_{C}(S))
$$
and a traceless \emph{parabolic Higgs bundle} $(E_{\bf v},\theta)$ on $\big(C, S\big)$ consists of a quasiparabolic vector bundle $E_{\bf v}$
together with a traceless Higgs field $\theta$ on $E_{\bf v}$. It is $\mu_{\cA}$-\textit{semistable} (respectively $\mu_{\cA}$-\textit{stable}) if for every  line subbundle $L \subset E$ invariant under $\theta$, we have $\mu_{\cA}(L,E) \leq  \mu_{\cA}(E)$
(respectively $\mu_{\cA}(L,E) <  \mu_{\cA}(E)$).

Given $d\in \mathbb Z$ let $\mathcal{H}_{\cA}(d)$ be the moduli space parametrizing  ${\rm S}$-equivalence classes of traceless $\mu_{\cA}$-semistable parabolic Higgs bundles  $(E_{\bf v},\theta)$ on $\big(C, S\big)$ with $\deg(E)=d$.  For a fixed line bundle $L$ of degree $\deg(L)=d$, we denote by $\mathcal{H}_{\cA}(L)$ the subvariety of  $\mathcal{H}_{\cA}(d)$ given by those parabolic Higgs bundles  with fixed determinant $\det E=L$. The existence of these moduli spaces follows from \cite{Yo93, Yo95}. 

Similarly to parabolic vector bundles,  given $w\in C$ we assume that $L=\mathcal O_C(w)$ and write $\mathcal{H}_{\cA}$ for the corresponding moduli space 
\[
\mathcal{H}_{\cA} = \mathcal{H}_{\cA}(\mathcal O_C(w))
\]
and $\mathcal{H}_{\cA}^s $ for the locus of $\mu_{\cA}$-stable Higgs bundles. 
%


The moduli space $\mathcal{H}_{\cA}$ is a normal quasiprojective variety of dimension $2m$, where $m=\dim \mathcal{M}_{\cA}$. When $E_{\bf v}\in \mathcal{M}_{\cA}^s$ then it follows by \cite[Theorem 2.4]{Yo95} that 
$$
T_{E_{\bf v}}^*\mathcal{M}_{\cA}^s \cong \Gamma(\mathcal{SE}nd_0(E_{\bf v})\otimes \omega_{C}(S)).
$$

%
We can perform an elementary transformation $elem_{I,L_0}:\mathcal{H}_{\cA} \to\mathcal{H}_{\cA'}$, as described in Section \ref{quasiparvec}, on the pair $(E_{\bf v},\theta)$. For instance, since $\theta$  is strongly parabolic then ${\rm Res}(\theta, p_i)$ is nilpotent with 
respect to the parabolic direction $V_i$ of $E$, $i=1,\dots, n$, then its restriction to $E' \subset E$ induces a homomorphism 
\[
\theta': E' \to E'\otimes \omega_C(S)
\]
which is nilpotent with respect to the parabolic direction ${\bf v'}$ of $E'$. If $e_1$, $e_2$ are local sections which generate $E$ near a parabolic point $p_i$ with $e_1(p_i)\in V_i$ then $\theta$ is given by
\begin{eqnarray*}
\theta=\left(
\begin{array}{ccc} 
xa & b  \\
xc & -xa  \\
\end{array}
\right)
\end{eqnarray*}
while that  
\begin{eqnarray*}
\theta'=\left(
\begin{array}{ccc} 
xa & xb  \\
c & -xa  \\
\end{array}
\right).
\end{eqnarray*}

\end{say}

\begin{say}{\bf The Hitchin map and spectral curves.}\label{Hitchinspec}\rm
\;Given  $(E_{\bf v},\theta)\in \mathcal H_{\cA}$, 
since ${\rm Res}(\theta, p_i)$ is nilpotent for every parabolic point $p_i\in C$, then  $\det(\theta)$ lies in the linear subspace 
$$\Gamma(\omega_{C}^{\otimes 2}(S))\subset \Gamma(\omega_{C}^{\otimes 2}(2S))$$ consisting of sections of $\Gamma(\omega_{C}^{\otimes 2}(2S))$ vanishing at $p_1, \dots, p_n$.  The \textit{Hitchin map} is defined as
$$
\begin{array}{cccc}
\frak{h}: &\mathcal{H}_{\cA}& \longrightarrow & \Gamma(\omega_{C}^{\otimes 2}(S))\\
      & (E_{\bf v},\theta) & \longmapsto & \det(\theta).
\end{array}
$$



We shall describe the fibers of the Hitchin map, for this purpose we recall the definition of the spectral curve.  Denote by $P$ the total space of the sheaf $\omega_{C}(S)$, with natural map
${\bf q}: P\longrightarrow C$. 
Let 
$s\in {\rm H}^0(P, {\bf q}^*(\omega_{C}(S)))$
be the tautological section: $s(v)$ is $v$ itself. Given $a\in{\rm H}^{0}(C, \omega_{C}^{\otimes 2}(2S))$ and $\alpha\in{\rm H}^{0}(C, \omega_{C}(S))$, we define the {\it spectral curve} $C_t\subset P$, $t=(\alpha, a)$, as the zero locus  of the section 
\[
s^2 -{\bf q}^*(\alpha) \cdot s + {\bf q}^*(a) \in{\rm H}^0(P, {\bf q}^*(\omega_{C}^{\otimes 2}(2S))). 
\]
We note that $C_t$  is singular if and only if there is a multiple zero of $a$ which is a zero of $\alpha$.  In addition, it comes with a degree two map  
\[
q_t = {\bf q}|_{C_t}:C_t \longrightarrow C.
\]


An equivalent definition is as follows. We can define a  structure of commutative ring  on  $\mathcal O_C\oplus \omega_C(S)^*$ induced by $t = (\alpha, a)$:
\stepcounter{thm}
\begin{eqnarray}\label{strucring}
(a_0,a_1)\cdot (b_0,b_1) := (a_0b_0-a(a_1b_1),a_0b_1+a_1b_0-\alpha (a_1b_1)).
\end{eqnarray}
This makes $\mathcal O_C\oplus \omega_C(S)^*$ an $\mathcal O_{C}$-algebra, which will be denoted by $\mathcal A_t$, and it is locally given by 
\[
\mathcal A_t(U) \simeq \frac{\mathcal O_C(U)[z]}{(z^2 -\alpha z + a)}.
\] 
The spectral curve can be  defined as
$$
C_t = {\rm Spec}(\mathcal A_t).
$$ 
In particular we see that 
\stepcounter{thm}
\begin{eqnarray}\label{qsO}
(q_t)_*\mathcal O_{C_t} = \mathcal O_C\oplus \omega_C(S)^*.
\end{eqnarray}


The genus $g_t$ of $C_t$ is given by 
\stepcounter{thm}
\begin{eqnarray}\label{genusXs}
g_t =  4(g-1)+n+1.
\end{eqnarray}
See \cite[Remark 3.2]{BNR89}. 

%
%

From now on we assume $\alpha = 0 $ and write simply $C_a$.

\begin{Remark}\label{singelip}
For $a\in\Gamma(\omega_{C}^{\otimes 2}(S))$,  $C_a$ is singular if and only if either $a\in\Gamma( \omega_{C}^{\otimes 2}(S-p_i))$ or $a\in\Gamma( \omega_{C}^{\otimes 2}(S-2p))$ with $p\in C$  distinct of $p_i$.  When $C$ is elliptic, then $C_a$  is a smooth curve of genus $n+1$ for general $a\in\Gamma(\omega_{C}^{\otimes 2}(S))$. 
\end{Remark} 

\end{say}

\begin{say}{\bf The locus of singular spectral curves.}\label{discriminant}\rm
\;In this section we study the locus of characteristic polynomials with singular spectral curves. The pair $(C,S)$ can be recovered from this locus,  this is the approach of \cite[Section 5]{HR04}, \cite[Section 4]{BGM13} and \cite[Section 4]{AG19}.  

We shall use the following notation: the dual space $(\mathbb P^N)^{\vee}$ parametrizes the set of hyperplanes in $\mathbb P^N$ and given a subvariety $X$ of $\mathbb P^N$, its dual variety $X^{\vee}\subset (\mathbb P^N)^{\vee}$ is the closure of the set of hyperplanes containing the projective tangent space of a smooth point of  $X$. 

As before  $C$ denotes an elliptic curve and $S=p_1+\cdots+p_n$ an effective reduced divisor on it. 
 Following the previous section,  we denote by $H_i=\Gamma( \omega_{C}^{\otimes 2}(S-p_i))$ and $D_p=\Gamma(\omega_{C}^{\otimes 2}(S-2p))$. We know  that $C_a$ is singular if and only if $a\in \mathcal D$, where
$$\displaystyle{\mathcal D= \cup_{i}H_i \cup_{p\in C-S} D_p\subset \Gamma( \omega_{C}^{\otimes 2}(S))}.$$
See Remark \ref{singelip}. We call $\mathcal D$ the {\it locus of singular spectral curves}. Notice that  $D_{p_i}\subset H_i$ for every $p_i\in S$, and then we can write $\mathcal D$ as 
\[
\mathcal D = \cup_{i} H_i \cup \tilde{\mathcal D} 
\]
where $\tilde{\mathcal D}$ is the divisor formed by sections having a zero of multiplicity at least two
\[
\tilde{\mathcal D}  = \cup_{p\in C} D_{p}. 
\]

If $S$ has  degree $\deg(S)\ge 3$, then $\omega_C^{\otimes 2}(S)$ is very ample, i.e., it defines an embedding 
\begin{eqnarray*}
C \longrightarrow  \mathbb P\Gamma( \omega_{C}^{\otimes 2}(S))^{\vee}
\end{eqnarray*}
sending $p\in C$ to the point at  $\mathbb P\Gamma(\omega_{C}^{\otimes 2}(S))^{\vee}$ which corresponds to the hyperplane $\Gamma( \omega_{C}^{\otimes 2}(S-p))$. Therefore we identify $C$ with its image in $\mathbb P\Gamma( \omega_{C}^{\otimes 2}(S))^{\vee}$. Observe that $C^{\vee}=\P\tilde{\mathcal D}\subset \P \Gamma( \omega_{C}^{\otimes 2}(S))$, and therefore by biduality we conclude that
$$
(\P\tilde{\mathcal D})^{\vee} = C.
$$
Also we have $(\P {H_i})^{\vee} = p_i \in C$ and  then we recover  $(C, S)$ from $\mathcal D$. We summarize the above discussion in the following proposition.

\begin{Proposition}\label{dual}
Let $C$ be an elliptic curve and  assume $\deg(S)\ge 3$. Let $\mathcal D$ be locus of singular spectral curves and let $a\in \Gamma(\omega_{C}^{\otimes 2}(S))$.  Then:
\begin{enumerate}
\item The spectral curve $C_a$ is singular if and only if  $a\in \mathcal D$. If  $C_a$ is non-integral then $a\in\Gamma(\omega_C^{\otimes 2})$. 
 \item The line bundle $\omega_C^{\otimes 2}(S)$  defines an embedding 
$C \hookrightarrow \mathbb P\Gamma( \omega_{C}^{\otimes 2}(S))^{\vee}$. 
 \item\label{dualtorelli} There is a decomposition of $\mathcal D$ on irreducible components $\mathcal D = \cup_{i=1}^n H_i \cup \tilde{\mathcal D}$ such that  $(\P\tilde{\mathcal D})^{\vee} = C$ and  $(\P H_i)^{\vee} = p_i$. 
\end{enumerate}
\end{Proposition}

\end{say}

\begin{say}{\bf Fibers of the  Hitchin map.}\label{fibers}\rm
\;Let us  assume that  $C_a$ is irreducible and smooth.  If  $M$ is a line bundle on $C_a$  then $E=(q_a)_*(M)$ is a rank two vector bundle on $C$ and the tautological section $s_a=s_{|C_a}$ induces a homomorphism ${(s_a)}_{*}: E\longrightarrow E\otimes \omega_{C}(S)$.  We can compute the degree of $E=(q_a)_*M$ using (\ref{qsO}) and the following identity 
$$
{\rm det} (E) \simeq {\rm det}((q_a)_*\mathcal O_{C_a} )\otimes {\rm Nm}(M)
$$
where ${\rm Nm}(M)$ is the norm map (see \cite[Chap. IV Ex. 2.6]{Ha97}):
$$
\deg M = \deg (E) + n +2g - 2.
$$
Let us fix ${\bf d} = 1+n + 2g - 2$ and let ${\rm Pic}^{\bf d}(C_a)$ be the variety parametrizing isomorphism classes of line bundles on $C_a$ of degree ${\bf d}$. Using the correspondence that associates $M$ to the pair $(E, {(s_a)}_{*})$ one can show that the general fiber of the Hitchin map $\frak{h}$ is an Abelian variety, this is essentially a consequence of \cite[Proposition 3.6]{BNR89}. See also \cite{LM10}. 


\begin{Proposition}\label{general fiber}
Let $C$ be an elliptic curve and assume $\deg S=n\ge 3$.  Let $\mathcal H_{\cA_F}$  be the moduli space of  $\mu_{\cA_F}$-semistable traceless parabolic Higgs bundles on $(C,S)$ satisfying $\det(E)=\mathcal O_C(w)$.
\begin{enumerate}
\item\label{1} If  $a\in  \Gamma( \omega_{C}^{\otimes 2}(S))\backslash \mathcal D$, then the spectral curve $C_a$ is irreducible and smooth of genus $n+1$, the fiber $\frak{h}^{-1}(a)$ is contained in $\mathcal H_{\cA_F}^s$ and it is isomorphic to 
\[
{\rm Prym}(C_a/C) = \left\{ M\in {\rm Pic}^{n+1}(C_a)\;:\;\; \det((q_a)_*M)=\mathcal O_C(w)\right\}. 
\]
\item\label{2}  For  general $a\in \mathcal D$, $C_a$  is an integral curve whose only singularity is one simple node,   $\frak{h}^{-1}(a)$ is contained in $\mathcal H_{\cA_F}^s$ and it is a uniruled variety. 
\end{enumerate}
\end{Proposition}

\proof
Let us prove (\ref{1}). By Proposition \ref{dual}, the spectral curve is smooth and integral for $a\in  \Gamma( \omega_{C}^{\otimes 2}(S))\backslash \mathcal D$.  Thus from \cite[Proposition 3.6]{BNR89},  there is a bijective correspondence between  ${\rm Pic}^{n+1}(C_a)$  and isomorphism classes of pairs $(E, \theta)$ where $E$ is a vector bundle of rank two and degree $1$, and $\theta: E\to E\otimes \omega_C(S)$ is a homomorphism with $\tr(\theta)=0$ and $\det (\theta) = a$. 
We note that since $\det(\theta)$ lies in the subspace $\Gamma( \omega_{C}^{\otimes 2}(S))$ of $\Gamma( \omega_{C}^{\otimes 2}(2S))$ then any residual matrix ${\rm Res}(\theta,p_i)$ has determinant equal to zero. So, the parabolic direction $V_i$ is defined as the kernel of ${\rm Res}(\theta,p_i)$.
 We claim that $(E_{\bf v}, \theta)$ is $\mu_{\cA_F}$-stable. If not, then there is a line subbundle $L\subset E$ invariant under $\theta$ and  this implies that $C_a$ is non-integral, which gives a contradiction. To conclude the proof of $(\ref{1})$ we note that $\frak{h}^{-1}(a)$  is formed by those $(E_{\bf v}, \theta)$ with $\det(E)=\mathcal O_C(w)$.  


We now prove (\ref{2}).  We first note that   $C_a$ has a simple node at the point over $p_i$ for general  $a\in H_i=\Gamma( \omega_{C}^{\otimes 2}(S-p_i))$. Indeed,  since $\Gamma( \omega_{C}^{\otimes 2}(S-2p_i))$ has dimension $n-2$, then it is a proper subspace of  $H_i$. In addition, $C_a$ is singular over $p_j$, $j\neq i$, if and only if $a\in\Gamma( \omega_{C}^{\otimes 2}(S-p_i-p_j))$ and it is singular over a nonparabolic point $p$ if and only if $a\in\Gamma( \omega_{C}^{\otimes 2}(S-p_i-2p))$. 
We conclude that for general $a\in H_i$,  $C_a$  is an integral curve having only one singularity which is a simple node over $p_i$. The same reasoning can be applied to a general point of $\tilde{\mathcal D}$. Now let $a\in \mathcal D$ be general, as before we can see that  $\frak{h}^{-1}(a)$  is contained in $\mathcal H_{\cA_F}^s$ and  one obtains a bijective correspondence between $\frak{h}^{-1}(a)$ and the isomorphism classes of torsion free sheaves $M$ of rank $1$ and degree $n+1$ on $C_a$ such that $\det((q_a)_*M) = \mathcal O_C(w)$. Equivalently,  $\frak{h}^{-1}(a)$ is a fiber of the map $\psi:{\rm \overline{P}ic}^{n+1}(C_a)\to {\rm Pic }^1(C)$ that sends  $M$ to $\det((q_a)_*M)$. We now note that  $\frak{h}^{-1}(a)$ is uniruled. Indeed, by \cite[Proposition 2.2]{Bho92} the compactified Jacobian ${\rm \overline{P}ic}^{n+1}(C_a)$ is uniruled and  $\psi$ contracts any rational curve because there is no nonconstant morphism from a rational curve to an Abelian variety. 



%
\endproof

\end{say}

\section{{The Torelli theorem for the model {$\mathcal M_{\cA_F}$}}} \label{main_section} 

In this section we  prove Theorem \ref{main}, our main result. For this, we need some preliminary results. 

We let $E_1$ denote the unique, up to isomorphism,  nontrivial extension 
\[
0\to \mathcal O_C \to E_1 \to \mathcal O_C(w)\to 0.
\]
%
%


\begin{say}{\bf Indecomposable  parabolic vector bundles.}\rm\label{sub:indeunstable}
A quasiparabolic bundle $(E,\textbf{v})$ is decomposable if there exists a decomposition $E=L\oplus M$ such that each parabolic direction is contained either in $L$ or in $M$. In this case, we write 
$$
(E,\textbf{v}) = (L,{\bf{v_1}}) \oplus  (M,{\bf{v_2}}).
$$

\begin{Proposition}\label{classi_ind_uns}
Let $C$ be an elliptic curve and assume $\deg S=n\ge 2$.  Let  $E_{\bf v}$ be an  indecomposable parabolic bundle which is not $\mu_{\cA_F}$-stable. Then there exists an elementary transformation which transforms $E_{\bf v}$ into one of the following:
\begin{enumerate}
	\item\label{case1ind}  $E=E_1$ and for any  destabilizing subbundle $L\subset E_1$,  the number of parabolic directions away from $L$ is at most $\frac{n}{2} -1$;
	\item\label{case2dec} $E=L\oplus L^{-1}(w)$, $L^2\simeq \mathcal O_C(w+p_i)$ for some $i=1, \dots, n$, $V_i$ does not lie in $L$, there is no embedding of $L^{-1}(w)$  containing $V_i$ and all the other parabolic directions lie in $L$. 
\end{enumerate}
\end{Proposition}

\proof 
First let us assume  $E=E_1$. 
Since  $E_{\bf v}$ is not $\mu_{\cA_F}$-stable,   there exists  a   line subbundle $L\subset E_1$ satisfying 
\begin{eqnarray}\label{m0m1}
1-2\deg(L) +\frac{m_0}{2} - \frac{m_1}{2} \le 0
\end{eqnarray}
where $m_0$ is the number of parabolic directions not lying in $L$ and $m_1$ is the number of parabolic directions that lie in $L$. Maximal subbundles of $E_1$ have degree zero, then $\deg L \le 0$ and from (\ref{m0m1}) we get $m_0\le\frac{n}{2} -1$.

We now assume that $E_{\bf v}$ is indecomposable as parabolic  bundle, but $E=L\oplus L^{-1}(w)$. Since  $L\oplus L^{-1}(w)\simeq M\oplus M^{-1}(w)$ with $M = L^{-1}(w)$ we may assume $\deg L =k \ge 1$. Thus $L$ is the only maximal subbundle and it corresponds to a section of $\mathbb PE$ with self-intersection $1-2k$. In order to arrive in cases (\ref{case1ind}) or (\ref{case2dec}) of the statement,  we will perform an elementary transformation $elem_I$  over parabolic points $\{p_i\}_{i\in I}$ whose parabolic directions are away from $L$ and such that $I\subset \{1,\dots, n\}$ has even cardinality.

First we will show that we can assume $\deg L=1$, up to elementary transformation.  The family of embeddings $L^{-1}(w)\hookrightarrow L\oplus L^{-1}(w)$ is parametrized by $\Gamma( L^{2}(-w))$, thus from Riemann-Roch theorem we conclude that it has dimension $2k-1$. Assuming  $k>1$,   given a parabolic direction $V_i$ outside $L$  we can choose an embedding of $L^{-1}(w)$ passing through it, because  $\Gamma( L^{2}(-w))$ has no base points. In the same way given a set $\{V_{i_1}, \dots, V_{i_{2k-2}}\}$ of $2k-2$ parabolic directions outside $L$ we can choose an embedding of $L^{-1}(w)$ passing through them. In particular, since $E_{\bf v}$ is indecomposable there are at least $2k-1$ parabolic directions,  $\{V_{i_1}, \dots, V_{i_{2k-2}}, V_{i_{2k-1}}\}$ away from $L$.  Choosing  $I=\{i_1,\dots, i_{2k-2}\}$, the elementary transformation $elem_I$ transforms $L$ to a subbundle whose corresponding section has self-intersection $-1$ and the transformed of  $L^{-1}(w)$ corresponds to a section of self-intersection $+1$. Thus after elementary transformation we can assume $E=L\oplus L^{-1}(w)$ with $\deg L = 1$.  

When $E=L\oplus L^{-1}(w)$ with $\deg L = 1$, the family of embeddings $L^{-1}(w)\hookrightarrow L\oplus L^{-1}(w)$ is one dimensional. This family corresponds to a one dimensional family of sections which have self-intersection $+1$. Writing $L^2(-w) = \mathcal O_C(p)$ we see that the linear system $\Gamma( L^{2}(-w))$ has a base point at $p$, which means that the family of $+1$ sections given by $L^{-1}(w)$ has a base point at the fiber over $p$.  Thus given a parabolic direction $V_i$ outside $L$, we can find an embedding of $L^{-1}(w)$ passing through it if and only if $p\neq p_i$. On the one hand, since $n\ge 2$, if $p\neq p_i$ then there exists at least one parabolic $V_{i}\subset L^{-1}(w)$ and, since $E_{\bf v}$ is indecomposable, at least one parabolic $V_j \nsubseteq L^{-1}(w)$,  $V_j \nsubseteq L$.  In this setting we perform one elementary transformation $elem_{i,j}$ in order to transform $E$ into $E_1$. On the other hand, if $p=p_i$ either we are in case (\ref{case2dec}) of the statement or we can transform $E$ into $E_1$ applying an elementary  transformation $elem_{i,j}$ over two parabolic directions outside $L$.  This finishes the proof of the proposition.
  
\endproof

\end{say}

\begin{say}{\bf The complement of $T^*\mathcal M_{\cA_F}^s$.}\rm\label{complement big 2}
In this section we study the complement of $T^*\mathcal M_{\cA_F}^s$ in $\mathcal H_{\cA_F}^s$. We will need the following result which concerns to Higgs bundles whose underlying parabolic bundle is decomposable.

\begin{Lemma}\label{casedecomposable}
Let $C$ be an elliptic curve  and assume $\deg S\ge 2$. Let $(E_{\bf{v}}, \theta)$ be a $\mu_{\cA_F}$-semistable traceless Higgs bundle whose  underlying parabolic bundle $(E,\bf{v})$ is decomposable  $(E,{\bf{v}}) = (L,{\bf{v_1}}) \oplus  (L^{-1}(w),{\bf{v_2}})$, $\deg L\ge 1$. Let $D$ be the divisor corresponding to those parabolic directions  which lie in $L$ and let $m=\deg D$. Then 
\[
2\deg L\le n-m+1
\] 
and we have:
\begin{enumerate}
 \item if $L^2\simeq \mathcal O_C(S-D+w)$ then  $\Gamma(\mathcal{SE}nd_0(E_{\bf v})\otimes \omega_{C}(S))$ has dimension $n+2$;
 \item otherwise,  $\Gamma(\mathcal{SE}nd_0(E_{\bf v})\otimes \omega_{C}(S))$ has dimension $n+1$. 
\end{enumerate}
\end{Lemma}

\proof

Let us  assume that $E$ is defined by the cocycle 
\begin{eqnarray*}
G_{ij}=\left(
\begin{array}{ccc} 
g_{ij} & 0  \\
0 & f_{ij}  \\
\end{array}
\right)
\end{eqnarray*}
where $\{g_{ij}\}$ determines the line bundle $L$ and $\{f_{ij}\}$ defines $L^{-1}(w)$. 
A traceless  Higgs field $\theta$ on $E$  with logarithmic poles at $S$ is determined by  $\theta_i$ in charts $U_i\subset C$ where 
\begin{eqnarray*}
\theta_i=\left(
\begin{array}{ccc} 
\alpha_i & \beta_i  \\
\gamma_i & -\alpha_i  \\
\end{array}
\right)
\in {\rm GL}_2(\omega_{U_i}(S))
\end{eqnarray*}
with $\{\alpha_i\}$, $\{\beta_i\}$ and $\{\gamma_i\}$ satisfying the compatibility conditions
$$
\theta_i\cdot G_{ij}=G_{ij}\cdot \theta_j 
$$
on each intersection $U_i\cap U_j$. Equivalently, 
\begin{displaymath}
\left\{ \begin{array}{ll}
\alpha = \{\alpha_i\} \in \Gamma(\omega_C(S))\\
\{\beta_i\}  \;\; \text{induces an element}\;\; \beta = \{f_ig_i^{-1}\beta_i\}\in  \Gamma(\omega_C(S)\otimes L^{2}(-w)) \\
\{\gamma_i \} \;\; \text{induces an element}\;\; \gamma = \{g_if_i^{-1}\gamma_i\}\in \Gamma(\omega_C(S)\otimes L^{-2}(w)).
\end{array} \right.
\end{displaymath}
where $f_{ij}=f_i/f_j$ and $g_{ij}=g_i/g_j$ are meromorphic resolutions of the cocycles.

Since any parabolic direction lies either in $L$ or in $L^{-1}(w)$, the condition which says that $\theta$ is strongly parabolic yields
\begin{displaymath}
\left\{ \begin{array}{ll}
\alpha \in \Gamma(\omega_C)\\
\beta \in  \Gamma(\omega_C(D)\otimes L^{2}(-w)) \\
\gamma \in \Gamma(\omega_C(S-D)\otimes L^{-2}(w)).
\end{array} \right.
\end{displaymath}

If $\deg L > (n-m+1)/2$ then $\Gamma(\omega_C(S-D)\otimes L^{-2}(w)) = \{0\}$. This implies $\gamma=0$ and  $L$ invariant under $\theta$. Since we have $\mu_{\cA_F}(L,E) >  \mu_{\cA_F}(E)$ one concludes that  $(E_{\bf v}, \theta)$ is not $\mu_{\cA_F}$-semistable and we get a contradiction. This finishes the first assertion of the statement. 

An element $\theta$ of  $\Gamma(\mathcal{SE}nd_0(E_{\bf v})\otimes \omega_{C}(S))$ is determined by $\alpha, \beta$ and $\gamma$ as above, thus  the conclusion of the proof follows from Riemann-Roch theorem, which can be applied to $\Gamma(\omega_C)$, $\Gamma(\omega_C(D)\otimes L^{2}(-w))$ and $ \Gamma(\omega_C(S-D)\otimes L^{-2}(w))$. 


\endproof

The next result shows that under hypothesis of Lemma \ref{casedecomposable} we can assume  $(E,{\bf{v}}) = (L,{\bf{v_1}}) \oplus  (L^{-1}(w),{\bf{v_2}})$ with $\deg L = 1$, up to an elementary transformation. 

\begin{Lemma}\label{decdeg1}
Let $C$ be an elliptic curve  and assume $\deg S=n\ge 2$. Let $(E_{\bf{v}}, \theta)$ be a $\mu_{\cA_F}$-semistable traceless Higgs bundle whose  underlying parabolic bundle $E_{\bf{v}}$ is decomposable $(E,{\bf{v}}) = (L,{\bf{v_1}}) \oplus  (L^{-1}(w),{\bf{v_2}})$, $\deg L\ge 1$. Then up to an elementary transformation we can assume $\deg L=1$. 
\end{Lemma}

\proof
By Lemma \ref{casedecomposable}, we obtain $m\le n+1-2\deg L$. Then there are at least $2\deg L -1$ parabolic directions away from $L$, i.e., lying   in  $L^{-1}(w)$.  After performing an elementary transformation over $2\deg L - 2$ of them, the transformed of $L$ is a line subbundle of degree $1$.
\endproof

\begin{Lemma}\label{caseindecomposable}
Assume we are in case (\ref{case2dec}) of Proposition  \ref{classi_ind_uns}, then any element $\theta\in\Gamma(\mathcal{SE}nd_0(E_{\bf v})\otimes \omega_{C}(S))$ leaves $L$ invariant.  In particular, $\theta$ is not $\mu_{\cA_F}$-semistable. 
\end{Lemma}

\proof
 For simplicity, we assume $i=1$, i.e., $E=L\oplus L^{-1}(w)$ and $L^2\simeq \mathcal O_C(w+p_1)$. As in the proof of Lemma \ref{casedecomposable}, a traceless  Higgs field $\theta$ on $E=L\oplus L^{-1}(w)$  with logarithmic poles at $S$ is determined by  a family of $\{\theta_i\}$ 
\begin{eqnarray*}
\theta_i=\left(
\begin{array}{ccc} 
\alpha_i & \beta_i  \\
\gamma_i & -\alpha_i  \\
\end{array}
\right)
\in {\rm GL}_2(\omega_{U_i}(S))
\end{eqnarray*}
satisfying the compatibility conditions $\theta_i\cdot G_{ij}=G_{ij}\cdot \theta_j $. 
%

Since any parabolic direction lies in $L$, unless $V_1$, the condition of being strongly parabolic at $p_j$, $j\neq 1$, together with $L^2\simeq \mathcal O_C(w+p_1)$  yield
\begin{displaymath}
\left\{ \begin{array}{ll}
\alpha \in \Gamma(\omega_C(p_1))\\
\beta=\{f_ig_i^{-1}\beta_i\} \in  \Gamma(\omega_C(S+p_1)) \\
\gamma=\{g_if_i^{-1}\gamma_i\} \in \Gamma(\omega_C).
\end{array} \right.
\end{displaymath}
Now since $\Gamma(\omega_C(p_1)) = \Gamma(\omega_C)$ the residual matrix of $\theta$ at $p_1$ is given by 
\begin{eqnarray*}
{\rm Res}(\theta; p_1)=\left(
\begin{array}{ccc} 
0 & b \\
c & 0 \\
\end{array}
\right)
\end{eqnarray*}
Finally, since the parabolic direction $V_1$ does not lie in $L$ neither in $L^{-1}(w)$ we get $b=c=0$. This implies $\beta=\{f_ig_i^{-1}\beta_i\} \in  \Gamma(\omega_C(S)) $ and $\gamma=\{g_if_i^{-1}\gamma_i\} \in \Gamma(\omega_C(-p_1))=\{0\}$. In particular, $\gamma = 0$ and this yields  $L$ invariant under $\theta$. 

\endproof

\begin{Proposition}\label{compTM}
Let $C$ be an elliptic curve and assume $\deg S=n\ge 3$.  The complement of $T^*\mathcal M_{\cA_F}^s$ in $\mathcal H_{\cA_F}^s$ has codimension at least $\frac{n}{2}$. 
\end{Proposition}

\proof
Let $Z=\mathcal H_{\cA_F}^s\setminus T^*\mathcal M_{\cA_F}^s$ be the complement. Elements of $Z$ correspond to those traceless $\mu_{\cA_F}$-stable parabolic Higgs bundles $(E_{\bf v}, \theta)$ whose underlying parabolic vector bundle  $E_{\bf v}$ is  not $\mu_{\cA_F}$-stable.  The finite group ${\bf El}$ formed by elementary transformations, described in Section \ref{quasiparvec}, acts on $\mathcal H_{\cA_F}$ preserving the locus $Z$. Let us consider the following families of semistable parabolic Higgs bundles: 
\begin{enumerate}
\item $Z_{dec} = \{ (E_{\bf v}, \theta)\in Z\;:\; E_{\bf v} \;\text{is decomposable}\;E_{\bf{v}} = (L,{\bf{v_1}}) \oplus  (L^{-1}(w),{\bf{v_2}}), \deg L =1\}.$
\item $Z_{ind} = \{(E_{\bf v}, \theta)\in Z\;:\;  E=E_1  \}.$
\end{enumerate} 

Proposition \ref{classi_ind_uns} classifies  indecomposable parabolic bundles which are not  $\mu_{\cA_F}$-stable, up to the action of ${\bf El}$. In addition, by Lemma \ref{caseindecomposable}  we can exclude case (\ref{case2dec}) of Proposition \ref{classi_ind_uns}. When $E_{\bf v}$ is decomposable $(E,{\bf{v}}) = (L,{\bf{v_1}}) \oplus  (L^{-1}(w),{\bf{v_2}})$, it follows from Lemma \ref{decdeg1} that we may assume $\deg L=1$, up to an elementary transformation.  Therefore, 
we conclude that 
\[
Z = \bigcup_{\varphi\in {\bf El}}\varphi(Z_{dec})\cup \varphi(Z_{ind}).
\]



When  $(E_{\bf v}, \theta)\in Z_{dec}$ 
the underlying  parabolic bundle   $E_{\bf{v}} = (L,{\bf{v_1}}) \oplus  (L^{-1}(w),{\bf{v_2}})$ is completely determined by $(L,{\bf{v_1}})$,  $\deg L= 1$, because $L$ is the only maximal subbundle. Let  $D$ be the divisor  formed by parabolic points whose parabolic directions  belong to ${\bf{v_1}}$ and let $m=\deg D$ ($m=0$ included). In order to bound the dimension of $Z_{dec}$ we may assume that $D$ is fixed.  The space of Higgs bundles over $E_{\bf{v}}$ is given by the quotient of  $\Gamma(\mathcal{SE}nd_0(E_{\bf v})\otimes \omega_{C}(S))$ by the group of automorphisms of $E$ fixing parabolic directions, namely the group  
\begin{eqnarray*}
{\rm Aut}(E_{\bf v}) = \left\{
\sigma_{a,c} = \left(
\begin{array}{ccc} 
a & 0  \\
0 & c  \\
\end{array}
\right)
\;:\;\; a,c\in\mathbb C^*
\right\}.
\end{eqnarray*}
By Lemma \ref{casedecomposable},   $\Gamma(\mathcal{SE}nd_0(E_{\bf v})\otimes \omega_{C}(S))$ has dimension $n+1$, unless  $L^2\simeq \mathcal O_C(S-D+w)$ which gives dimension $n+2$.  
The group ${\rm Aut}(E_{\bf v})$ acts on $\Gamma(\mathcal{SE}nd_0(E_{\bf v})\otimes \omega_{C}(S))$ with positive dimensional generic orbit. For instance, $\sigma_{a,c}$  acts as a dilatation on $\beta$ and $\gamma$,  sending $\beta$ to $(ac^{-1})\beta$ and $\gamma$ to $(a^{-1}c)\gamma$. Hence, for each $L\in {\rm Pic}^1(C)\simeq C$  the quotient of $\Gamma(\mathcal{SE}nd_0(E_{\bf v})\otimes \omega_{C}(S))$ by ${\rm Aut}(E_{\bf v})$ has dimension at most $n$, unless $L$ is a square root of $ \mathcal O_C(S-D+w)$ which gives dimension $n+1$. We conclude that $Z_{dec}$ has dimension at most $n+1$. 

We now assume  that $(E_{\bf v}, \theta)\in Z_{ind}$, i.e., we are in case (\ref{case1ind}) of Proposition \ref{classi_ind_uns}. 
Since $E_1$ has no automorphisms, besides trivial ones,  the parabolic vector bundle $(E_1,{\bf v})$  is completely determined by a destabilizing subbundle $L$ of degree $k\le 0$ and by a set $\{V_{i_1}, \dots, V_{i_{m_0}}\}$ of $m_0$ parabolic directions outside $L$, $0\le m_0\le\frac{n}{2}-1$.  Let us fix $k\le 0$ and a subset  $I=\{i_1,\dots, i_{m_0}\}$ of  $\{1, \dots, n\}$. Let $I^c\subset \{1, \dots, n\}$ be the complement of $I$ and let $Z_{ind}(k,I^c)$ be the family of parabolic Higgs bundles $(E_1,{\bf v}, \theta)$ having a destabilizing subbundle $L$ of degree $k$ satisfying  $V_j\subset L$ for all $j\in I^c$.  For each $L\in   {\rm Pic}^k(C)$, the underlying  parabolic bundle $(E_1,{\bf v})$ in this family is determined by those parabolic directions over $I$: 
\[
(V_{i_1}, \dots, V_{i_{m_0}}) \in \mathbb P (E_{i_1})\times \cdots \times \mathbb P (E_{i_{m_0}})\simeq (\mathbb P^1)^{m_0}. 
\] 
Varying $L\in   {\rm Pic}^k(C)$, these parabolic bundles form a family of dimension $1+m_0$. We note that each $(E_1,{\bf v})$ is $\mu_{\bf a }$-stable with respect to a weight ${\bf a} = (a_1,\dots ,a_n)$ satisfying $0<a_i<1$ and  $\sum a_i<1$. This gives 
\[
\Gamma(\mathcal{SE}nd_0(E_1,{\bf v})\otimes \omega_{C}(S))\simeq T_{(E_1,{\bf v})}^* \mathcal M_{\bf a} 
\]
and hence $\Gamma(\mathcal{SE}nd_0(E_1,{\bf v})\otimes \omega_{C}(S))$ has dimension $n$. This shows that $Z_{ind}(k,I^c)$ has dimension at most $1+m_0+n$, with $m_0\le\frac{n}{2}-1$. Therefore, since  $Z_{ind}$ is a countable union of $Z_{ind}(k,I^c)$,  we obtain  that its codimension is at least $\frac{n}{2}$. This concludes the proof of the proposition. 

%
\endproof

\end{say}

\begin{say}\rm{\bf Affinization.} 
Let $X$ be an algebraic variety over $\mathbb C$ and $\mathbb C[X]$ its algebra of global regular functions. There is a map $X\to {\rm Spec} (\mathbb C[X])$, called {\it affinization map}, which sends $x\in X$ to the maximal ideal $\frak{m}_x\subset\mathbb C[X]$ formed by global regular functions vanishing at $x$. If $U\subset X$ is an affine open subset then the restriction of the affinization map to $U$ is the map $U \to {\rm Spec}\mathbb C[X]$ induced by the restriction homomorphism 
\[
r_U: \mathbb C[X] \to \mathbb C[U]
\] 
which sends $f$ to  $f|_U$. 

We will show that  the ring $A=\C[T^*\mathcal M_{\cA_F}^s]$  produces a map
\[
\tilde{h}: T^*\mathcal M_{\cA_F}^s \longrightarrow {\rm Spec}(A)
\]
which turns out to be the Hitchin map, up to an automorphism on the basis:

\begin{Proposition}\label{affinization}
Let $C$ be an elliptic curve and assume $\deg S\ge 3$. We have an isomorphism ${\rm Spec}(A)\simeq \Gamma(\omega_{C}^{\otimes 2}(S))$ and the Hitchin map $h: T^*\mathcal M_{\cA_F}^s\longrightarrow \Gamma(\omega_{C}^{\otimes 2}(S))$ is simply the affinization 
\[
\tilde{h}: T^*\mathcal M_{\cA_F}^s \longrightarrow {\rm Spec}(A)
\]
up to this isomorphism. 
\end{Proposition}
\proof
Let $\Sigma  = \Gamma(\omega_{C}^{\otimes 2}(S))$. Given a regular function $f\in \C[\Sigma]$,  the composition with $h$ gives a regular function  $f\circ h\in A$. Since $\Sigma={\rm Spec}(\C[\Sigma])$, this gives a map from  ${\rm Spec}(A)$ to $\Sigma$. Reciprocally, given $g\in A$, it follows from Proposition \ref{compTM} and Hartogs theorem that $g$ extends to a regular function $\overline{g}$ on $\mathcal H_{\cA_F}^s$. By Proposition \ref{general fiber}, the restriction of $\overline{g}$ to a fiber $\frak{h}^{-1}(a)$, of the Hitchin map $\frak{h}: \mathcal H_{\cA_F}^s \to \Sigma$, is constant,  varying $a$ away from a codimension two subset of $\Sigma$.  
%
Then there is a regular function $f:\Sigma \longrightarrow \C$ satisfying $g = f\circ h$.  This shows that the isomorphism ${\rm Spec}(A)\simeq \Sigma$  is induced by the homomorphism $h^*:\mathbb C[\Sigma] \to A$ which sends $f$ to $f\circ h$. 

Now if $U\subset T^*\mathcal M_{\cA_F}^s$ is an affine open subset, since $\tilde{h}|_U$ is induced by the restriction homomorphism $r_U: A \to \mathbb C[U]$,  the conclusion of the proof of the proposition follows from the commutativity  of the diagram
\[
 \xymatrix{
 \mathbb C[U]  \ar@{<-}[d]_{r_U}  \ar@{<-}[dr]^{(h|_U)^*}   \\
   A \ar@{<-}[r]_{h^*} & \mathbb C[\Sigma].}
\]

\endproof

\end{say}

\begin{say}{\bf Torelli theorem.}\rm 
 \;We now prove the Torelli type theorem, which is  Theorem \ref{main} of the introduction: 

\begin{thm}\label{main inside}
Let $C$ and $C'$ be two elliptic curves. Let $S=p_1+\cdots+p_n$ and $S'=p'_1+\cdots+p'_n$ denote  two effective reduced divisors on $C$ and $C'$, respectively, with $n\ge 3$. Given $w\in C$ and $w'\in C'$, consider the two moduli spaces $\mathcal M_{\cA_F}$ and $\mathcal M'_{\cA_F}$ of $\mu_{\cA_F}$-semistable parabolic vector bundles,  over $(C,S)$ and $(C',S')$, with fixed determinant line bundle $\mathcal O_C(w)$ and $\mathcal O_{C'}(w')$,  respectively. If these moduli spaces are isomorphic, then there is an isomorphism $C\longrightarrow C'$ sending $S$ to $S'$. 
\end{thm}

\proof
Let $\phi: \mathcal M_{\cA_F} \longrightarrow \mathcal M'_{\cA_F}$ be an isomorphism and consider the induced isomorphism on the cotangent bundles $\psi:T^*\mathcal M_{\cA_F}^s \longrightarrow T^*{\mathcal M'}_{\cA_F}^s$. By Proposition \ref{affinization}, there is a morphism of affine varieties $\xi:{\Gamma(\omega_C^{\otimes 2}(S))}\longrightarrow{\Gamma(\omega_C^{\otimes 2}(S'))}$ making the following diagram commute:
  \[
  \begin{tikzpicture}[xscale=4.5,yscale=-1.5]
    \node (A0_0) at (0, 0) {$T^{*}{\mathcal M}_{\cA_F}^s$};
    \node (A0_1) at (1, 0) {$T^{*}{\mathcal M'}_{\cA_F}^s$};
    \node (A1_0) at (0, 1) {${\Gamma(\omega_C^{\otimes 2}(S))}$};
    \node (A1_1) at (1, 1) {${\Gamma(\omega_C^{\otimes 2}(S'))}$.};
    \path (A0_0) edge [->]node [auto] {$\scriptstyle{\psi}$} (A0_1);
    \path (A1_0) edge [->]node [auto] {$\scriptstyle{\xi}$} (A1_1);
    \path (A0_1) edge [->]node [auto] {$\scriptstyle{h'}$} (A1_1);
    \path (A0_0) edge [->,swap]node [auto] {$\scriptstyle{h}$} (A1_0);
  \end{tikzpicture}
  \]
We note that the $\C^*$-action by dilatations on the fibers of the map $T^{*}\mathcal{M}_{\cA_F}^s \longrightarrow \mathcal{M}_{\cA_F}^s$ induces a $\C^*$-action on ${\Gamma(\omega_C^{\otimes 2}(S))}$: $(\lambda, a) \mapsto \lambda^2\cdot a$. Since $\psi$ is $\C^*$-equivariant then  $\xi$ is also $\C^*$-equivariant, and then linear. 
By Proposition \ref{general fiber}, the fiber $h^{-1}(a)$ is an open subset of an abelian variety for  $a\in \Gamma(\omega_C^{\otimes 2}(S))\backslash\mathcal D$ and it is an open subset of a uniruled variety for general $a\in\mathcal D$. The same assertion holds for $h'$.  Since abelian varieties are not uniruled, one obtains  that $\xi$ preserves the locus of singular spectral curves, {\it i.e.}, it sends $\mathcal D$ to $\mathcal D'$.  Therefore, it induces a map $\xi^{\vee}: \P{\Gamma(\omega_C^{\otimes 2}(S))}^{\vee}\longrightarrow \P{\Gamma(\omega_{C'}^{\otimes 2}(S'))}^{\vee}$ that sends $\P{\mathcal D}^{\vee}$ to $\P{\mathcal D'}^{\vee}$. By Proposition \ref{dual}(\ref{dualtorelli}), $\xi^{\vee}$ sends $(C,S)$ to $(C',S')$. 
\endproof

\end{say}

\bibliographystyle{amsalpha}
\bibliography{Biblio}
\end{document}